\documentclass[12pt]{article}
\usepackage{amsmath}
\newcommand{\be}{\begin{equation}}
\newcommand{\ee}{\end{equation}}
\newcommand{\ba}{\begin{eqnarray}}
\newcommand{\ea}{\end{eqnarray}}
\newcommand{\ban}{\begin{eqnarray*}}
\newcommand{\ean}{\end{eqnarray*}}

 \newcommand{\qed}{\hspace*{\fill}\rule{3mm}{3mm}\quad}
\newcommand{\Pf}{\noindent  {\em Proof.} }

\newcommand{\sect}[1]{\section{#1}  \setcounter{equation}{0}}

\newtheorem{lem}{Lemma}[section]
 
\begin{document}
\newtheorem{defn}[lem]{Definition}
\newtheorem{theo}[lem]{Theorem}
\newtheorem{cor}[lem]{Corollary}
\newtheorem{prop}[lem]{Proposition}
\newtheorem{rk}[lem]{Remark}
\newtheorem{ex}[lem]{Example}
\newtheorem{note}[lem]{Note}
\newtheorem{conj}[lem]{Conjecture}

\title{The Logarithmic Sobolev Inequality Along The Ricci Flow: The Case $\lambda_0(g_0)=0$}
\author{Rugang Ye \\ {\small Department  of Mathematics} \\
{\small University of California, Santa  Barbara}}
\date{August 1, 2007}
\maketitle

\sect{Introduction}

In [Y1] and [Y2], logarithmic Sobolev inequalities along the Ricci flow 
in all dimensions $n \ge 2$ were obtained using Perelman's entropy monotonicity, which lead to Sobolev inequalities and $\kappa$-noncollpasing estimates. In particular, a uniform logarithmic Sobolev inequality, a uniform 
Sobolev inequality and a uniform $\kappa$-noncollapsing estimate were obtained without any restriction on 
time, provided that the smallest eigenvalue $\lambda_0(g_0)$ of the operator $-\Delta+\frac{R}{4}$ for the initial metric  
is positive.  In this paper, we extend these uniform results to the case $\lambda_0(g_0)=0$. 
 
Consider a compact manifold $M$ of dimension $n \ge 2$.  Let $g=g(t)$ be a smooth solution of the Ricci flow 
\ba
\frac{\partial g}{\partial t}=-2Ric
\ea
on $M \times [0, T)$ for some (finite or infinite) $T>0$ with a given initial metric 
$g(0)=g_0$.\\

\noindent {\bf Theorem A}  {\it Assume that $\lambda_0(g_0)=0$. 
For each $t\in [0, T)$ and each $\sigma>0$ there holds 
\ba \label{sobolevA}
\int_M u^2 \ln u^2 dvol \le \sigma \int_M (|\nabla u|^2 +\frac{R}{4}u^2) dvol
-\frac{n}{2}\ln \sigma +C
\ea
for all $u \in W^{1,2}(M)$ with $\int_M u^2 dvol=1$, where $C$ depends only on $(M, g_0)$. 
} \\

Note that as in [Y1] a log gradient version of the logarithmic Sobolev inequality follows as a consequence of 
(\ref{sobolevA}). We omit its statement. Our next result provides a dependence of the above $C$ on $g_0$ in terms of rudimentary geometric data. 
To simplifiy the statements, we assume that $|Rm| \le 1$ for $g_0$, which can 
always be achieved by a rescaling. Then we can also assume $T \ge 2\alpha(n)$ for a positive 
constant $\alpha(n)$ depending only on $n$ such that $|Rm| \le 2$ on $[0, \alpha(n)]$. (Namely the maximal possible $T$ such that the solution 
$g=g(t)$ can be extended to a smooth solution of the Ricci flow on $[0, T)$ has this property.)
\\

\noindent {\bf Theorem B}  {\it  There are for each 
$v_0>0$, each $D_0>0$, each $\epsilon>0$ and each integer $l \ge 3$ a positive number  $C=C(v_0, D_0, \epsilon, l,  n)$ 
with the following properties. Assume $\lambda_0(g_0)=0, vol_{g_0}(M) \ge v_0, diam_{g_0} \le D_0$ and the normalization conditions $|Rm|_{g_0} \le 1$ and 
$T \ge 2 \alpha(n)$. $vol_{g_0}(M) \ge v_0$. Then one of the following two cases must occur: \\
1) $g(\alpha(n))$ lies in 
the $\epsilon$-neighborhood of a Ricci flat metric on $M$ in the $C^l$ norm,\\
2) the logarithmic Sobolev inequality (\ref{sobolevA})  holds 
true for each $t \in [0, T)$, each $\sigma>0$ and all $u \in W^{1,2}(M)$ with 
$\int_M u^2 dvol=1$, where $C=C(v_0, D_0, \epsilon, l, n)$.  }\\

It turns out that we have a better result in dimension $n=3$. The same holds true in 
dimension $n=2$. But Theorem E and Theorem 3.7 in [Y2] provide a stronger result in this 
dimension.  \\

\noindent {\bf Theorem C} {\it Assume that $n=3$.   There is for each 
$v_0>0$ and  each $D_0>0$   a positive 
number $C=C(v_0, D_0)$  with the following properties. Assume $\lambda_0(g_0)=0, vol_{g_0}(M) \ge v_0, diam_{g_0} \le D_0$ and the normalization conditions $|Rm|_{g_0} \le 1$ and 
$T \ge 2 \alpha(n)$.  Then 
the logarithmic Sobolev inequality (\ref{sobolevA})  holds 
true for each $t \in [0, T)$, each $\sigma>0$ and all $u \in W^{1,2}(M)$ with 
$\int_M u^2 dvol=1$, where $C=C(v_0, D_0)$.  }\\

As in [Y1], Theorem A, Theorem B and Theorem C lead to Sobolev inequalities along the Ricci flow, which in turn lead to 
$\kappa$-noncollpasing estimates. We consider only the case $n \ge 3$ although the methods also work for 
$n=2$, because the case 
$n=2$ is covered by the results in [Y2].   \\

\noindent {\bf Theorem D} {\it Assume that $n \ge 3$ and $\lambda_0(g_0)=0$. Then there holds for each 
$t \in [0, T)$
\ba
\label{sobolevD}
\left( \int_M |u|^{\frac{2n}{n-2}} dvol \right)^{\frac{n-2}{n}} \le 
A\int_M (|\nabla u|^2+\frac{R}{4}u^2)dvol+B \int_M u^2 dvol
\ea
for all $u \in W^{1,2}(M)$, where $A$ and $B$ depend only on $(M, g_0)$. 
We have $A=A(v_0, D_0)$ and $B=B(v_0, D_0)$ for given $v_0>0$ and $ D_0>0$, if $n=3$ and 
$g_0$ satisfies the conditions in Theorem C.  We also have $A=A(v_0, D_0, n)$ and $B=B(v_0, 
D_0, n)$ for given $v_0>0$ and $ D_0>0$, provided that
$g_0$ satisfies  the conditions in Theorem B and $g(\alpha(n))$ does not lie in 
the $\epsilon$-neighborhood of any Ricci flat metric on $M$ in the $C^3$ norm. }
\\

\noindent {\bf Theorem E} {\it Assume that $n=3$ and $\lambda_0(g_0)=0$.   Let $L>0$ and $t \in [0, T)$. Consider the Riemannian manifold 
$(M, g)$ with $g=g(t)$. Assume $R\le \frac{1}{r^2}$ on a geodesic ball $B(x, r)$ with $0<r \le L$. Then 
there holds 
\ba \label{noncollapse}
vol(B(x, r)) \ge \left(\frac{1}{2^{n+3}A+2BL^2}\right)^{\frac{n}{2}} r^n,
\ea
where $A$ and $B$ are from Theorem D. }\\

As in [Y1], the above results extend to various versions of the modified Ricci flows. Moreover, 
the $\kappa$-noncollapsing estimates ensure that we can obtain smooth blow-up limits at the 
time infinity under the assumption that $\lambda_0(g_0)=0$.  We omit the statements of those results because they 
are completely analogous to the corresponding ones in [Y1].

\sect{The Proofs}

\noindent {\bf Proof of Theorem A}  Consider a fixed $t_1 \in (0, T)$. Let $u_1$ be a positive eigenfunction for the eigenvalue $\lambda_0(g(t_1))$ 
associated with the metric $g(t_1)$, such that 
$\int_M u_1^2 dvol=1$ with respect to $g(t_1)$.  Let $f=f(t)$ be the smooth solution of the 
equation 
\ba \label{f}
\frac{\partial f}{\partial t}=-\Delta f+|\nabla f|^2-R
\ea
on $[0, t_1]$ with $f(t_1)=-2\ln u_1$.  Note that (\ref{f}) is equivalent to
\ba
\frac{\partial v}{\partial t}=-\Delta v+Rv,
\ea
where $v=e^{-f}$. So the solution $f(t)$ exists. We also infer 
$\frac{d}{dt}\int_M v dvol=0$, and hence $\int_M v dvol =1$ for all $t \in [0, t_1]$. 
 
 We set $u=e^{-\frac{f}{2}}.$  By [P, (1.4)] we then have
\ba
\frac{d}{dt} \int_M (|\nabla u|^2+\frac{R}{4}u^2)dvol 
=\frac{1}{4} \frac{d}{dt}\int_M (|\nabla f|^2+R)e^{-f}dvol \ge \frac{1}{2} \int_M |Ric+\nabla^2 f|^2 e^{-f}dvol.
\nonumber \\
\ea
It follows that 
\ba \label{solitoninequality}
\lambda_0(g(t_1)) 
\ge \lambda_0(g_0) + \frac{1}{2} \int_0^{t_1} \int_M |Ric+\nabla^2 f|^2 e^{-f}dvol dt.
\ea

We choose $t_1=\min\{\frac{T}{2}, 1\}.$  If $\lambda_0(g(t_1)>0$, we first apply 
Theorem A in [Y1] or Theorem A in [Y2] to obtain (\ref{sobolevA}) for $g(t)$ on $[0, t_1]$.
Then we apply Theorem 4.2 in [Y1] to obtain (\ref{sobolevA}) on 
$[t_1, T)$ with a larger $C$.   If $\lambda_0(g(t_1))=0$, we deduce from 
(\ref{solitoninequality}) 
\ba
Ric+\nabla^2 f=0
\ea
on $[0, t_1]$. It follows that $g=g(t)$ is a steady Ricci soliton for all $t$. Hence the logarithmic 
Sobolev inequality for $g_0$ provided by Theorem 3.3 in [Y1] holds true for all $g(t)$. 
Actually, [CK, Proposition 5.20] implies that $g_0$ is Ricci flat, hence  
$g(t)=g_0$ for all $t$. 
\qed \\

\begin{lem} \label{mulemma} For given $v_0>0, D_0>0, \epsilon>0$ and $l \ge 3$ 
there is a positive constant $\mu_0=\mu_0(v_0, D_0, \epsilon, 
l, n)$ with the following properties. Let $g=g(t)$ be 
a smooth solution of the Ricci flow on $M \times [0, T)$ which satisfies the 
normalization conditions in Theorem B and the condtions $vol_{g(0)} \ge v_0$ and 
$diam_{g(0)} \le D_0$. Assume $\lambda_0(g_0)=0$.  If $\lambda(g(\alpha(n))) < \mu_0$, then 
$g(\alpha(n))$ lies in the $\epsilon$-neighborhood of a Ricci flat metric with 
respect to the $C^l$-norm. 
\end{lem}
\Pf Assume that $\mu_0$ does not exist. Then we can find a sequence of 
manifolds $M_k$ of a fixed dimension $n$ and a sequence of smooth solutions $g_k=g_k(t)$ 
on $M_k \times [0, T_k)$ satisfying the normalization conditions and the conditions 
$vol_{g_k(0)} \ge v_0$ and $diam_{g_k(0)} \le D_0$, such that $\lambda_0(g_k(0))=0$,
$\lambda_0(g_k(\alpha(n))) \rightarrow 0$, and $g_k(\alpha(n))$ does not lie in 
the $\epsilon$-neighborhood of any Ricci flat metric with respect to the $C^l$-norm. 
By Gromov-Cheeger-Hamilton compactness theorem [H], we can find a subsequence of 
$(M_k, g_k)$, which we still denote by $(M_k, g_k)$, such that 
$(M_k, g_k, [\frac{\alpha(n)}{2}, \alpha(n)])$ converge smoothly to a limit Ricci flow 
$(M, g,  [\frac{\alpha(n)}{2}, \alpha(n)])$.  There holds $\lambda_0(g(\alpha(n)))=0$. 
By the monotonicity of $\lambda_0$ along the Ricci flow (see [P] or the above proof of 
Theorem A), we have $\lambda_0(g_k(t)) \ge 0$ for all $t \in [0, T_k)$. Hence 
$\lambda_0(g(t)) \ge 0$ for all $t \in [\frac{\alpha(n)}{2}, \alpha(n)]$.  
Now the argument in the proof of Theorem A implies that $g(\alpha(n))$ 
is Ricci flat. But $(M_k, g_k(\alpha(n)))$ converge smoothly to $(M, g(\alpha(n))$, so 
$g_k(\alpha(n))$ lies in the $\epsilon$-neighborhood of a Ricci flat metric 
with respect to  the $C^{l}$-norm whenever $k$ is large enough. This is 
a contradiction.  \qed \\

\noindent {\bf Proof of Theorem B}   Assume that $g(\alpha(n))$ does not lie in 
the $\epsilon$-neighborhood of any Ricci flat metric with respect to 
the $C^l$-norm. Then $\lambda_0(g(\alpha(n))) \ge \mu_0$ by 
Lemma \ref{mulemma}.    Now we obtain a desired logarithmic Sobolev inequality for 
$g(t)$ on $[0, \alpha(n)]$ by Theorem A in [Y1]. Alternatively, we can 
apply the arguments in [Y3] for controlling the evolution of the Sobolev constant to bound the Sobolev constant for $g(t)$ on 
$[0, \alpha(n)]$, and then apply Theorem 3.3 in [Y1] to infer the desired 
logarithmic Sobolev inequality.  Next we apply the bound for the Sobolev constant 
at $t=\alpha(n)$, the bound $\lambda_0 \ge \mu_0$ at $t=\alpha(n)$ and 
the arguments in the proof of Theorem 3.5 in [Y1] to deduce a logarithmic Sobolev inequality of the kind 
[Y1, (3.11)] at $t=\alpha(n)$.  Then we apply Theorem B in [Y1] on $[\alpha(n), T)$ 
and 
combine it with Theorem A in [Y1]. Then we arrive at the desired logarithmic 
Sobolev inequality on $[\alpha(n), T)$. \qed \\

\noindent {\bf Proof of Theorem C}  By [GIK], there is for each given flat metric $g$
an $\epsilon$-neighbord of $g$ with respect to the $C^6$-norm, such that 
the Ricci flow starting at any metric in the neighborhood converges 
smoothly to a Ricci flat metric at a fixed exponential rate as $t \rightarrow \infty$.
Moreover, the limit Ricci flat metric lies in the $2\epsilon$-neighborhood of $g$ with respect to 
the $C^6$-norm. 
We call such a neighborhood a {\it Ricci contraction  $\epsilon$-neighborhood}. 

Now for given $v_0>0, D_0>0$ and $K_0>0$ the moduli space ${\mathcal M}^0(v_0, 
D_0, K_0)$ of flat metrics on $M$ with $vol\ge v_0, diam \le D_0$ and 
$|Rm| \le K_0$  is $C^{\infty}$ compact modulo diffeomorphisms by Gromov-Cheeger compactness theorem and 
the Einstein equation.  So there is a uniform $\epsilon=\epsilon(v_0, D_0, K_0)$ such that 
each $g \in {\mathcal M}^0(v_0, D_0, K_0)$ has a Ricci contraction $\epsilon$-neighborhood. 
Moreover, there is a unform upper bound 
$C(v_0, D_0, K_0)$ for the Sobolev constant for $g \in {\mathcal M}^0(v_0, D_0, K_0)$.

Now consider for given $v_0>0$ and $D_0>0$ a smooth solution of the Ricci flow $g=g(t)$ 
satisfying the normalization conditions and the conditions $vol_{g_0}(M) \ge v_0$ and 
$diam_{g_0}(M) \le D_0$.  Let $\epsilon>0$ and $l=6$, where $\epsilon$ is to be determined.  Assume that $g(\alpha(n))$ lies 
in the $\epsilon$-neighborhood of a Ricci flat metric $\bar g$ with respect to the $C^6$ norm.
Since $n=3$, $\bar g$ is flat.
There is a positive number  $\epsilon_0=\epsilon_0(v_0, D_0)$ such that 
$vol_{\bar g}(M) \ge \frac{1}{2}v_0, diam_{\bar g}(M) \le 2 D_0$ and 
$|Rm|_{\bar g} \le 3$, i.e. 
$\bar g \in {\mathcal M}^0(\frac{1}{2}v_0, 2D_0, 3)$, whenever $\epsilon \le \epsilon_0$.  Now we choose $\epsilon=
\min\{\epsilon_0(v_0, D_0), \epsilon(\frac{1}{2}v_0, 2D_0, 4)\}$.  Then 
the maximally extended $g(t)$ converges smoothly to a flat metric $g^*$  at exponential rate where the rate depends only on 
$v_0$ and $D_0$.  We can choose $\epsilon_0$ sufficiently small such that 
$g^* \in {\mathcal M}^0(\frac{1}{3}v_0, 3D_0, 4)$. We can also make 
the $C^6$ norm of $g(t)-g^*$ sufficiently small for all $t \ge \alpha(n)$ such that 
the Sobolev constant of $g(t)$ is bounded above by $2C(\frac{1}{3}v_0, 3D_0, 4)$ for all
$t \ge \alpha(n)$. A desired logarithmic Sobolev inquality then follows for 
$t \in [\alpha(n), T)$.   

A desired logarithmic Sobolev inequality for $g(t)$ on $[0, \alpha(n)]$ follows from Theorem A in 
[Y1]. It also follows from the arguments for controlling the evolution of the Sobolev constant 
in [Y3].  \qed \\

\end{document}